\documentclass[12pt,a4paper]{article}
\usepackage{multicol,amsmath,graphicx,amsfonts,amssymb,amsthm,tabularx}
\usepackage{bm,color,latexsym,mathrsfs,fancybox,framed,enumerate}
\usepackage{bbm}
\newcommand{\ind}[1]{\mathbbm{1}_{\{#1\}}}
\allowdisplaybreaks[1]
\numberwithin{equation}{section}
\theoremstyle{plain}
\newtheorem{Thm}{Theorem}[section]

\newtheorem{Lem}[Thm]{Lemma}
\newtheorem{Prop}[Thm]{Proposition}

\newtheorem{Rem}[Thm]{Remark}
\newcommand{\B}{\mathbb{B}}
\newcommand{\Z}{\mathbb{Z}}

\newcommand{\R}{\mathbb{R}}
\newcommand{\Bd}{\B^d}
\newcommand{\Zd}{\Z^d}
\newcommand{\E}{\mathbb{E}}
\newcommand{\nn}{\nonumber}
\newcommand{\Prob}{\mathbb{P}}
\newcommand{\Xvec}{\boldsymbol{X}}
\newcommand{\Yvec}{\boldsymbol{Y}}
\newcommand{\saw}[1]{\Omega(#1)}
\newcommand{\sss}{\scriptscriptstyle}
\newcommand{\gSAW}{\hat\Omega^{\sss\sf{good}}_{\delta,\Xvec}}

\newcommand{\vep}{\varepsilon}
\definecolor{dgreen}{rgb}{0,0.5,0}

\newcommand{\ha}[1]{h^\mathsf{a}_{#1}}
\newcommand{\hq}[1]{\hat h^\mathsf{q}_{#1}}
\newcommand{\Proof}[1]{\paragraph{\it #1}}
\newcommand{\QED}{\hspace*{\fill}\rule{7pt}{7pt}\smallskip}
\title{The quenched critical point for self-avoiding walk on random conductors}
\author{
Yuki~Chino\footnote{{\tt chino@math.sci.hokudai.ac.jp}}~~ and~
Akira~Sakai\footnote{{\tt sakai@math.sci.hokudai.ac.jp}}\\
Department of Mathematics\\
Hokkaido University
}
\begin{document}
\maketitle
\begin{abstract}
Following similar analysis to that in Lacoin \cite{lacoin2014non}, we can show
that the quenched critical point for self-avoiding walk on random conductors on
$\Zd$ is almost surely a constant, which does not depend on the location of the
reference point.  We provide upper and lower bounds which are valid for all
$d\ge1$.
\end{abstract}
\section{Introduction}
\subsection{Background}
Self-avoiding walk (SAW) is a statistical-mechanical model for chain-like
solvents and linear polymers.  SAW was first introduced by
Flory~\cite{flory1949conf,flory1953principles} in order to model and
investigate the behavior of polymer chains.  Since then,
many rigorous mathematical results on SAW have been proven,
while physicists have much more conjectures that are believed to be true.
Most of them are supported by numerical simulations and physical ideas
that have not been fully justified mathematically.

It would be more natural to consider an inhomogeneous environment in which
polymers lie.  In recent years, various models of SAW in a quenched random
environment have attracted attention of chemists, physicists and
mathematicians
\cite{chakra1981the,chakra1984static,harris1983self,meir1989self}.
One of them is SAW on a randomly diluted lattice, introduced
by Chakrabarti and Kert\'{e}sz~\cite{chakra1981the}.
Le Doussal and Machta~\cite{doussal1991self} investigated it by applying
a renormalization method on a hierarchical lattice and came up to some
conjectures.  Lacoin~\cite{lacoin2014existence} answered affirmatively to
one of them by showing that, on an infinite supercritical percolation cluster
in 2 dimensions, the quenched critical point
(defined by divergence of the quenched susceptibility) is strictly smaller
than the annealed one (defined by divergence of the average susceptibility).

In this paper, we investigate SAW in a different type of random environment,
which is topologically regular, but random in energy landscape.
The goal is to achieve better understanding of how the introduction of
randomness changes the properties of the critical point.

\subsection{The model and the main theorem}
Let $\Omega(x,y)$ be the set of (nearest-neighbor) self-avoiding paths
on $\Zd$ from $x$ to $y$, and let $\Omega(x)=\bigcup_{y\in\Zd}\Omega(x,y)$.
Denoting the length of $\omega$ by $|\omega|$ (i.e., $|\omega|=n$ for
$\omega=(\omega_0,\dots,\omega_n)$) and the energy cost of a bond between
consecutive monomers by $h\in\R$, we define the susceptibility as
\begin{align}
\chi_h=\sum_{\omega\in\Omega(x)}e^{-h|\omega|},
\end{align}
which is independent of the location of the reference point $x\in\Zd$.
Two other key observables are the number of $n$-step SAWs and
the two-point function:
\begin{align}\label{Ghcn}
c(n)=\sum_{\omega\in\Omega(x)}\ind{|\omega|=n},&&
G_h(x)=\sum_{\omega\in\Omega(o,x)}e^{-h|\omega|},
\end{align}
where $o$ is the origin of $\Zd$ and $\ind{\cdots}$ is the indicator function.
Obviously,
\begin{align}
\chi_h=\sum_{n=0}^\infty e^{-hn}c(n)=\sum_{x\in\Zd}G_h(x).
\end{align}
Due to subadditivity of $\log c(n)$, we can readily show that
$\chi_h<\infty$ if and only if $h>\log\mu$, where $\mu$ is the connective
constant for SAW \cite{madras2013self}:
\begin{align}\label{mu}
\mu=\lim_{n\to\infty}c(n)^{1/n}=\inf_nc(n)^{1/n}.
\end{align}
Therefore, $h_0\equiv\log\mu$ is the critical point of the susceptibility.
Many rigorous results on the behavior of these observables around
$h=h_0$ have been proven, especially in high dimensions $d>4$,
with the help of the lace expansion \cite{brydges1985self,madras2013self}.
However, there still remain many challenging open problems in two and three
dimensions.  See \cite{slade2011self} and the references therein.

Next, we introduce randomness to the environment.  Let $\Bd$ denote the set
of nearest-neighbor bonds in $\Zd$, and let $\Xvec=\{X_b\}_{b\in\Bd}$ be a
collection of integrable random variables whose law $\Prob$ is
translation-invariant and ergodic.
From a physical point of view, $X_b$ can be regarded as the magnitude of
resistance of a conductor attached to $b\in\Bd$, and therefore it is more
natural to assume $X_b\ge0$.  However, the results in this paper are all valid
without this assumption.  Given the environment $\Xvec$ and the strength of
randomness $\beta\ge0$, we define the quenched susceptibility at $x\in\Zd$ as
\begin{align}
\hat\chi_{h,\beta,\Xvec}(x)=\sum_{\omega\in\Omega(x)}e^{-\sum_{j
 =1}^{|\omega|}(h+\beta X_{b_j})},
\end{align}
where
\begin{align}
b_j\equiv b_j(\omega)=(\omega_{j-1},\omega_j).
\end{align}
Because of the inhomogeneity of $\Xvec$, the quenched susceptibility is no
longer translation invariant and does depend on the location of the reference
point $x$.  Similarly to the homogeneous case, we also define
\begin{align}
\hat c_{\beta,\Xvec}(x;n)&=\sum_{\omega\in\Omega(x)}
 e^{-\beta\sum_{j=1}^{|\omega|}X_{b_j}}\ind{|\omega|=n},\\
\hat G_{h,\beta,\Xvec}(x,y)&=\sum_{\omega\in\Omega(x,y)}
 e^{-\sum_{j=1}^{|\omega|}(h+\beta X_{b_j})}.
\end{align}
These quantities are reduced to $\chi_h$, $c(n)$ and $G_h(y-x)$, respectively,
when $\beta=0$.  Moreover,
\begin{align}
\hat\chi_{h,\beta,\Xvec}(x)=\sum_{n=0}^\infty e^{-hn}\hat c_{\beta,\Xvec}(x;n)
 =\sum_{y\in\Zd}\hat G_{h,\beta,\Xvec}(x,y).
\end{align}
Since $\hat\chi_{h,\beta,\Xvec}(x)$ is monotonic in $h$, we can define the
quenched version of the critical point as
\begin{align}
\hq{\beta,\Xvec}(x)=\inf\{h\in\R:\hat\chi_{h,\beta,\Xvec}(x)<\infty\}.
\end{align}

Our goal is to understand how the randomness of the environment $\Xvec$
affects the behavior of these quenched observables around the critical point.
There are numerous examples in which the introduction of randomness alters
the behavior of relevant observables.  Classical examples are Sinai's
one-dimensional random walk in a random medium \cite{sinai1983rwre} and Smith
and Wilkinson's branching processes in random environments \cite{sw1969bpre}.
More recent examples are the random pinning models
\cite{den2009random,giacomin2007random} and the directed polymer models
\cite{comets2004probabilistic}.

As a first step to understand the properties of the random variable
$\hq{\beta,\Xvec}(x)$, we consider the mean-field approximation
(or often called the annealing),
i.e., to take the average of $\hat\chi_{h,\beta,\Xvec}(x)$ over the
environment $\Xvec$.  Let
\begin{align}
\ha\beta=\{h\in\R:\E[\hat\chi_{h,\beta,\Xvec}(x)]<\infty\},
\end{align}
where $\E$ is the expectation for $\Prob$.  Since
$\Prob$ is translation-invariant, the annealed critical point
$\ha\beta$ does not depend on the location of the reference point
$x\in\Zd$.
We note that $\hq{\beta,\Xvec}(x)\le\ha\beta$ by definition.
In particular, if $\Xvec$ is i.i.d.~and the Laplace transform
\begin{align}
\lambda_\beta=\E[e^{-\beta X_b}]
\end{align}
exists, then
we can directly compute $\E[\hat c_{\beta,\Xvec}(x;n)]$ as
\begin{align}
\E[\hat c_{\beta,\Xvec}(x;n)]=\sum_{\omega\in\Omega(x):|\omega|=n}
 \,\prod_{j=1}^n\E[e^{-\beta X_{b_j}}]=\lambda_\beta^n\,c(n),
\end{align}
and the annealed susceptibility
$\E[\hat\chi_{h,\beta,\Xvec}(x)]$ as
\begin{align}\label{susca}
\E[\hat\chi_{h,\beta,\Xvec}(x)]=\sum_{n=0}^{\infty}e^{-hn}\E
 [\hat c_{\beta,\Xvec}(x;n)]&=\sum_{n=0}^{\infty}e^{-(h-\log\lambda_\beta)n}
 c(n)\nn\\
&=\chi_{h-\log\lambda_\beta}.
\end{align}
Therefore, 
\begin{align}\label{ha}
\ha\beta=h_0+\log\lambda_\beta.
\end{align}
By Jensen's inequality, we immediately see
\begin{align}\label{habd}
\ha{\beta}\ge h_0-\beta\E[X_b],
\end{align}
where the gap is $O(\beta^2)$ as $\beta\to0$.

The following theorem is the main result of this paper.


\begin{Thm}\label{thm:main}
Let $d\ge1$ and $\beta\ge0$.
The quenched critical point $\hq{\beta,\Xvec}(x)$ is
almost surely an $x$-independent constant.
Moreover, by abbreviating $\hq{\beta,\Xvec}(x)$ as $\hq{\beta,\Xvec}$,
we have
\begin{align}\label{hqbd}
h_0-\beta\E[X_b]\le\hq{\beta,\Xvec}\le\ha\beta,\quad\text{almost surely}.
\end{align}
For $d=1$, in particular, the lower bound is an equality.
\end{Thm}

Before proving this theorem in the next section, we give two remarks.

\begin{Rem}[On the first inequality in (\ref{hqbd})]
{\rm
For $d=1$, $\hq{\beta,\Xvec}=-\beta\E[X_b]$  (recall that $h_0=0$ for $d=1$)
is due to the fact that $c(n)$ is always two: either to the left or to the
right of the reference point.  Let
$h=-\beta\E[X_b]+\delta$ and
$\Delta_j=X_{(x+j-1,x+j)}-\E[X_b]$.  Then, we have
\begin{align}
\hat\chi_{h,\beta,\Xvec}(x)=1+\sum_{n=1}^\infty e^{-\delta n}\Big(e^{-\beta
 \sum_{j=1}^n\Delta_j}+e^{-\beta\sum_{j=0}^{n-1}\Delta_{-j}}\Big).
\end{align}
By applying the individual ergodic theorem to those two sequences
$\{\Delta_j\}_{j=1}^\infty$ and $\{\Delta_{-j}\}_{j=0}^\infty$, we can conclude
that the above series almost surely converges if and only if $\delta>0$.

For $d\ge2$, however, since $c(n)$ grows exponentially, it is hard to control
the speed of convergence along those walks at the same time.  Because of this
entropic effect, we strongly believe that the first inequality in (\ref{hqbd})
is a strict inequality.  So far, we have been able to prove it to be true only
for SAW on i.i.d.~random conductors in a homogeneous tree of degree $\ell\ge3$, if $\beta$ is
sufficiently small.  In fact, we can show the equality
$\hq{\beta,\Xvec}=\ha\beta$ as follows.  First, we set
$\delta=\ha\beta-h>0$ and use $\mu=\ell-1$ to obtain the rewrite
\begin{align}
\hat\chi_{h,\beta,\Xvec}(x)=\sum_{\omega\in\Omega(x)}e^{\delta|\omega|}
 \prod_{j=1}^{|\omega|}e^{-\ha\beta-\beta X_{b_j}}=1+\frac\ell{\ell-1}\sum_{n
 =1}^\infty e^{\delta n}Z_{\beta,\Xvec}(x;n),
\end{align}
where
\begin{align}
Z_{\beta,\Xvec}(x;n)=\sum_{\omega\in\Omega(x):|\omega|=n}\frac1{\ell(\ell
 -1)^{n-1}}\prod_{j=1}^n\frac{e^{-\beta X_{b_j}}}{\lambda_\beta}
\end{align}
is a positive martingale, and thus the limit
$Z_{\beta,\Xvec}(x;\infty)\equiv\lim_{n\to\infty}Z_{\beta,\Xvec}(x;n)$ exists
almost surely.  Adapting the statement in \cite[p.886]{johnson2011tree} to our
setting, we have the following dichotomy:
\begin{align}
f(\beta)\equiv\ha\beta-(\ha\beta)'\beta
 \begin{cases}
 >0\quad\Rightarrow\quad
  \Prob\big(Z_{\beta,\Xvec}(x;\infty)>0\big)=1,\\[5pt]
 \le0\quad\Rightarrow\quad
  \Prob\big(Z_{\beta,\Xvec}(x;\infty)=0\big)=1.
 \end{cases}
\end{align}
Since $f(0)=\log(\ell-1)>0$ and $f'(\beta)=-(\ha\beta)''\beta\le0$, we have
$f(\beta)>0$ for sufficiently small $\beta$, hence
$\hat\chi_{h,\beta,\Xvec}(x)=\infty$ almost surely.

Now we are back on $\Zd$.  If $\beta$ is large and $\E[X_b]>0$, then the gap
between the lower and upper bounds in (\ref{hqbd}) is large, and the inequality
(\ref{hqbd}) is no longer effective.  In the following specific case, however,
we may find a better bound.  Suppose that $\Prob(X_b=0)$ is bigger than
the critical point for oriented percolation on $\Zd_+$.  Then, there is almost
surely an $\Xvec$-free infinite oriented-percolation cluster $\mathcal{C}_x$
at some $x\in\Zd_+$, in which the number of $n$-step directed paths from $x$
grows exponentially in $n$ \cite[Theorem~3.1(2)]{fukushima2012on}.
The susceptibility $\hat\chi_{h,\beta,\Xvec}(x)$ can be bounded below by
restricting the sum over those directed paths in $\mathcal{C}_x$, implying
existence of a $\beta$-independent positive lower bound on $\hq{\beta,\Xvec}$.
}
\end{Rem}

\begin{Rem}[On the second inequality in (\ref{hqbd})]
{\rm
Although it is trivial by definition, the second inequality in (\ref{hqbd}) can
be proven in the following tedious way.  First, by the Markov inequality,
we have
\begin{align}\label{markov1}
\Prob\Big(\hat c_{\beta,\Xvec}(x;n)\ge n^2\E[\hat c_{\beta,
 \Xvec}(x;n)]\Big)\le\frac1{n^2}.
\end{align}
Then, by the Borel-Cantelli lemma, we can conclude that the opposite inequality
$\hat c_{\beta,\Xvec}(x;n)\le n^2\E[\hat c_{\beta,\Xvec}(x;n)]$
holds for all but finitely many $n$, implying almost sure convergence of
$\hat\chi_{h,\beta,\Xvec}(x)$ for $h>\ha\beta$.

We may improve it to a strict inequality in two dimensions by adapting the
idea of Lacoin \cite{lacoin2014non}.  In his setting (i.e., SAW on an infinite
supercritical percolation cluster in $\Z^2$), it is proven that there are
$b,\theta\in(0,1)$ such that
\begin{align}\label{framom}
\E[\hat c_{\beta,\Xvec}(x;n)^\theta]\le\big(b^n\E[\hat c_{\beta,
 \Xvec}(x;n)]\big)^\theta.
\end{align}
Then, by the Markov inequality, we have
\begin{align}\label{markov2}
\Prob\Big(\hat c_{\beta,\Xvec}(x;n)\ge n^{2/\theta}b^n\E
 [\hat c_{\beta,\Xvec}(x;n)]\Big)\le\frac1{n^2}.
\end{align}
By the Borel-Cantelli lemma again, we may conclude  $\hq{\beta}\le\ha{\beta}-\log\frac1b$.

Analyzing fractional moments, as in (\ref{framom}), has been a standard method
to investigate disordered systems. To see how it is used in other settings, we
refer to \cite{yilmaz2010differing} for random walks in random environments,
and to \cite{berger2010on,birkner2010annealed} for random pinning models.
}
\end{Rem}

\section{Proof of the main result}
We prove Theorem~\ref{thm:main} as follows.  In Section~\ref{ss:degeneration},
we prove the first half of Theorem~\ref{thm:main} by showing that the
quenched critical point is a degenerate random variable that does not depend on
the location of the reference point.
In Section~\ref{ss:lowerbd}, we complete the proof of Theorem~\ref{thm:main}
by showing the first inequality in (\ref{hqbd}).
Recall that its reduction to an equality for $d=1$ and the second
inequality in (\ref{hqbd}) for all $d\ge1$ have already been mentioned
in the previous section.

\subsection{Degeneration of the quenched critical point}\label{ss:degeneration}
Recall that $\Xvec=\{X_b\}_{b\in\Bd}$ is a collection of integrable (thus
almost surely finite) random variables whose law
$\Prob$ is translation-invariant and ergodic.  Following similar
analysis to that in Lacoin \cite{lacoin2014non}, we first prove that the
quenched critical point is independent of the location of  the reference point.


\begin{Lem}\label{lem:hq1}
The quenched critical point $\hq{\beta,\Xvec}(x)$ is almost surely a constant
function of $x\in\Zd$.
\end{Lem}


\Proof{Proof.}
We will show below that
\begin{align}\label{hqX}
\hat\chi_{h,\beta,\Xvec}(u)\le\hat\chi_{h,\beta,\Xvec}(v)^2+e^{h+\beta
 X_{(v,u)}}\hat\chi_{h,\beta,\Xvec}(v)
\end{align}
holds for any pair of neighboring vertices $u,v\in\Zd$.
Since $X_{(u,v)}$ is almost surely finite, it implies
$\hat\chi_{h,\beta,\Xvec}(u)<\infty$ if and only if
$\hat\chi_{h,\beta,\Xvec}(v)<\infty$.  Repeated applications of this inequality
to all neighboring vertices in $\Zd$, we conclude that all vertices are in the
same equivalent class, i.e., either $\hat\chi_{h,\beta,\Xvec}(x)<\infty$ for
all $x\in\Zd$ or $\hat\chi_{h,\beta,\Xvec}(x)=\infty$ for all $x\in\Zd$.
Therefore, $\hq{\beta,\Xvec}(x)$ does not depend on $x\in\Zd$,
almost surely.

It remains to show \eqref{hqX}.  First, we split the sum into two as
\begin{align}
\hat\chi_{h,\beta,\Xvec}(u)=\sum_{\omega\in\saw{u}}e^{-\sum_{j=1}^{|\omega|}
 (h+\beta X_{b_j})}\big( \ind{v\in\omega}+\ind{v\notin\omega}\big).
\end{align}
Due to subadditivity and reversibility, the contribution from
$\ind{v\in\omega}$ is bounded as
\begin{align}\label{up1}
\sum_{\omega\in\Omega(u):v\in\omega}e^{-\sum_{j=1}^{|\omega|}(h+\beta X_{b_j})}
&\le\underbrace{\sum_{\omega\in\saw{u,v}}e^{-\sum_{j=1}^{|\omega|}(h+\beta
 X_{b_j(\omega)})}}_{\hat G_{h,\beta,\Xvec}(u,v)}~\underbrace{\sum_{\eta\in
 \saw{v}}e^{-\sum_{j=1}^{|\eta|}(h+\beta X_{b_j(\eta)})}}_{\hat\chi_{h,\beta,
 \Xvec}(v)}\nn\\
&=\hat G_{h,\beta,\Xvec}(v,u)~\hat\chi_{h,\beta,\Xvec}(v)\nn\\[5pt]
&\le\hat\chi_{h,\beta,\Xvec}(v)^2.
\end{align}
On the other hand, by adding an extra step from $v$ to $u$, the contribution
from $\ind{v\notin\omega}$ is bounded as
\begin{align}\label{up2}
\sum_{\omega\in\saw{u}:v\notin\omega}e^{-\sum_{j=1}^{|\omega|}(h+\beta X_{b_j})}
&=e^{h+\beta X_{(v,u)}}\sum_{\omega\in\saw{u}:v\notin\omega}e^{-(h+\beta X_{(v,
 u)})}e^{-\sum_{j=1}^{|\omega|}(h+\beta X_{b_{j}(\omega)})}\nn\\
&=e^{h+\beta X_{(v,u)}}\sum_{\bar\omega\in\saw{v}:\bar\omega_1=u}
 e^{-\sum_{j=1}^{|\bar\omega|}(h+\beta X_{b_{j}(\bar{\omega})})}\nn\\
&\le e^{h+\beta X_{(v,u)}}\hat\chi_{h,\beta,\Xvec}(v),
\end{align}
where we have used the symmetry $X_{(u,v)}=X_{(v,u)}$ to form $\bar\omega$.
This completes the proof.
\QED


In the rest of this section, we simply denote $\hq{\beta,\Xvec}(x)$ by
$\hq{\beta,\Xvec}$.


\begin{Lem}\label{lem:hq2}
The quenched critical point $\hq{\beta,\Xvec}$ is
a degenerate random variable.
\end{Lem}


\Proof{Proof.}
Due to Lemma~\ref{lem:hq1}, the event $\{\hq{\beta,\Xvec}=h\}$ for any
$h\in\R$ is translation invariant.  Since $\Prob$ is ergodic, we can
conclude that $\Prob(\hq{\beta,\Xvec}=h)$ is either zero or one.
\QED

\subsection{Lower bound on the quenched critical point}\label{ss:lowerbd}
In this section, we prove the first inequality in (\ref{hqbd}) by showing
almost sure divergence of the quenched susceptibility
at $h=h_0-\beta\E[X_b]-\beta\delta$ for any $\beta>0$ and $\delta>0$.

Let $\Delta_b=X_b-\E[X_b]$ and define
\begin{align}
\Omega(x;n)&=\{\omega\in\Omega(x):|\omega|=n\},\\
\gSAW(x;n)&=\Big\{\omega\in\Omega(x;n):\big|\tfrac1n{\textstyle\sum_{j=1}^n}
 \Delta_{b_j(\omega)}\big|<\delta\Big\}.
\end{align}
Using this random set, we can bound $\hat\chi_{h,\beta,\Xvec}(x)$ at
$h=h_0-\beta\E[X_b]-\beta\delta$ as
\begin{align}\label{qsusc}
\hat\chi_{h,\beta,\Xvec}(x)=\sum_{\omega\in\Omega(x)}\frac1{\mu^{|\omega|}}
 e^{\beta|\omega|\big(\delta-\frac1{|\omega|}\sum_{j=1}^{|\omega|}\Delta_{b_j}
 \big)}
&\ge\sum_{n=1}^\infty\frac1{\mu^n}|\gSAW(x;n)|.
\end{align}
If there are infinitely many $n$ such that $|\gSAW(x;n)|\ge\frac12c(n)$,
then, by $c(n)\ge\mu^n$ (cf., (\ref{mu})), we obtain divergence of
the susceptibility.  Therefore,
\begin{align}\label{cond-on-good}
\Prob(\hat\chi_{h,\beta,\Xvec}=\infty)&\ge\underbrace{\Prob\Big(\hat\chi_{h,
 \beta,\Xvec}=\infty\Big|\limsup_{n\to\infty}\big\{|\gSAW(x;n)|\ge\tfrac12c(n)
 \big\}\Big)}_1\nn\\
&\quad\times\Prob\Big(\limsup_{n\to\infty}\big\{|\gSAW(x;n)|\ge\tfrac12c(n)
 \big\}\Big)\nn\\
&\ge\lim_{n\to\infty}\Prob\big(|\gSAW(x;n)|\ge\tfrac12c(n)\big).
\end{align}

To complete the proof, since $\Prob(\hat\chi_{h,\beta,\Xvec}(x)=\infty)$ is
either zero or one, it suffice to show that the rightmost limit is positive.
Here, we use the following Paley-Zygmund (PZ)
inequality~\cite{paley1933notes}: for a random variable $Z\ge0$
whose second moment is finite and for $\vep\in(0,1)$,
\begin{align}\label{PZineq}
\Prob(Z\ge\vep\E[Z])\ge(1-\vep)^2\frac{\E[Z]^2}{\E[Z^2]}.
\end{align}
Let $Z=|\gSAW(x;n)|$.  Notice that, by definition and ergodicity,
we can bound $\E\big[|\gSAW(x;n)|\big]$ from below as
\begin{align}\label{ergodicbd}
\E\big[|\gSAW(x;n)|\big]=\sum_{\omega\in\Omega(x;n)}\Prob\big(\big|{\textstyle
 \frac1n\sum_{j=1}^n}\Delta_{b_j(\omega)}\big|<\delta\big)\ge c(n)\big(1-o(1)
 \big).
\end{align}
Using this and the trivial inequality $\E\big[|\gSAW(x;n)|^2\big]\le c(n)^2$,
we obtain\footnote{One of two anonymous referees found the following much
simpler proof of (\ref{complementpz}).  First, by the trivial inequality
$|\gSAW(x;n)|\le c(n)$, we obtain
\begin{align}
\E\big[|\gSAW(x;n)|\big]
&\le\frac12c(n)\,\Prob\big(|\gSAW(x;n)|<\tfrac12c(n)\big)
 +c(n)\,\Prob\big(|\gSAW(x;n)|\ge\tfrac12c(n)\big)\nn\\
&=\frac12c(n)\Big(1+\Prob\big(|\gSAW(x;n)|\ge\tfrac12c(n)\big)\Big).
\end{align}
Combining this with (\ref{ergodicbd}), we can readily conclude
$\Prob\big(|\gSAW(x;n)|\ge\tfrac12c(n)\big)\ge1-o(1)$.
}
\begin{align}\label{complementpz}
\lim_{n\to\infty}\Prob\big(|\gSAW(x;n)|\ge\tfrac12c(n)\big)\ge\frac14>0,
\end{align}
as required.
\QED

\section{Another application of the PZ inequality}
Application of the PZ inequality is often dubbed the second-moment method.
It has been a standard tool to investigate disordered systems.  We show below
that the PZ inequality may also be used to investigate critical behavior for SAW on
i.i.d.~random conductors.  From now on, we assume that $\lambda_\beta<\infty$
for all $\beta\ge0$.

\begin{Prop}\label{prop:mfbd}
Suppose that
\begin{align}\label{bubble1}
B_1\equiv\E\bigg[\sum_{y\in\Zd}\hat G_{h,\beta,\Xvec}(x,y)^2\bigg]<\infty
\end{align}
and
\begin{align}\label{bubble2}
B_2\equiv\E\bigg[\sum_{y,z\in\Zd}\hat G_{h,\beta,\Xvec}(x,z)\,\hat G_{h,
 \beta,\Xvec}(z,y)^2\,\hat G_{h,\beta,\Xvec}(y,x)\bigg]<\infty
\end{align}
hold uniformly in $h>\ha\beta$.  Then, for any slowly-varying function
$L(h)\downarrow0$ as $h\downarrow\ha\beta$, we have
\begin{align}\label{highD}
\liminf_{h\downarrow\ha{\beta}}\Prob\bigg(\hat\chi_{h,\beta,\Xvec}(x)
 \ge\frac{L(h)}{h-\ha\beta}\bigg)\ge1-O(\beta^2).
\end{align}
\end{Prop}

Although the above result is conditional and still weak to establish a decisive
conclusion, it provides evidence to support the belief that, in high dimensions,
the coincidence $\hq{\beta,\Xvec}=\ha\beta$ occurs and the critical exponent
for $\hat\chi_{h,\beta,\Xvec}(x)$, if it exists, is bounded below by its
mean-field value 1.  For SAW in a homogeneous environment, the
conditions (\ref{bubble1})--(\ref{bubble2}) (in fact, the former implies the
latter because $B_2\le B_1^2$, which is a result of translation invariance and
the Cauchy-Schwarz inequality) are known to hold in dimensions $d>4$, via the
lace expansion \cite{brydges1985self,madras2013self}.  The lace expansion
yields a convolution equation for the two-point function, which is applicable
in both homogeneous and inhomogeneous settings.  In the current random
setting, however, because of the lack of translation invariance, we have not
been able to fully control the lace-expansion coefficients.  This is under
investigation in an ongoing project.

\Proof{Proof of Proposition~\ref{prop:mfbd}}
First, by replacing $Z$ in (\ref{PZineq}) by
$\hat\chi_{h,\beta,\Xvec}(x)$, we have
\begin{align}\label{disc}
\Prob\Big(\hat\chi_{h,\beta,\Xvec}(x)\ge\vep\E[\hat\chi_{h,
 \beta,\Xvec}(x)]\Big)\geq(1-\vep)^2\frac{\E[\hat\chi_{h,\beta,\Xvec}
 (x)]^2}{\E[\hat\chi_{h,\beta,\Xvec}(x)^2]}.
\end{align}
Since $\E[\hat\chi_{h,\beta,\Xvec}(x)]=\chi_{h-\log\lambda_\beta}$ (cf.,
(\ref{susca})) and $\chi_h\ge(h-h_0)^{-1}$ for all $h>h_0$ (cf., (\ref{mu})),
we have $\E[\hat\chi_{h,\beta,\Xvec}(x)]\ge(h-\ha\beta)^{-1}$ for all
$h>\ha\beta$.  Replacing $\vep$ in (\ref{disc}) by a slowly-varying function
$L(h)\downarrow0$ as $h\downarrow\ha\beta$, we can conclude (\ref{highD})
as soon as we can show
\begin{align}\label{2nd-1st^2:cond}
\frac{\E[\hat\chi_{h,\beta,\Xvec}(x)^2]-\E[\hat\chi_{h,\beta,
 \Xvec}(x)]^2}{\E[\hat\chi_{h,\beta,\Xvec}(x)]^2}\le O(\beta^2),
\end{align}
in the neighborhood of $\ha\beta$.

To prove (\ref{2nd-1st^2:cond}) under the assumptions
(\ref{bubble1})--(\ref{bubble2}), we introduce the notation
\begin{align}
H_{\Xvec}(\omega)=-\sum_{j=1}^{|\omega|}\big(h+\beta X_{b_j(\omega)}\big).
\end{align}
Let $\Yvec=\{Y_b\}_{b\in\Bd}$ be an independent copy of $\Xvec$.
Then, we obtain
\begin{align}\label{2nd-1st^2:eq1}
&\E[\hat\chi_{h,\beta,\Xvec}(x)^2]-\E[\hat\chi_{h,\beta,\Xvec}
 (x)]^2\nn\\
&=\sum_{\omega,\eta\in\Omega(x)}\E\bigg[e^{H_{\Xvec}(\omega)}
 \E_{\Yvec}\Big[e^{H_{\Xvec}(\eta)}-e^{H_{\Yvec}(\eta)}\Big]\bigg].
\end{align}
By the telescopic-sum representation, we can decompose
$e^{H_{\Xvec}(\eta)}-e^{H_{\Yvec}(\eta)}$ as
\begin{align}
e^{H_{\Xvec}(\eta)}-e^{H_{\Yvec}(\eta)}
&=\sum_{j=1}^{|\eta|}e^{H_{\Xvec}
 (\eta_{<j})}e^{-h}\Big(e^{-\beta X_{b_j(\eta)}}-e^{-\beta Y_{b_j(\eta)}}\Big)
 e^{H_{\Yvec}(\eta_{>j})},
\end{align}
where $\eta_{<j}=(\eta_0,\dots,\eta_{j-1})$ and
$\eta_{>j}=(\eta_{j+1},\dots,\eta_{|\eta|})$, with the convention
$H_{\Xvec}(\varnothing)=0$.  Substituting this back into (\ref{2nd-1st^2:eq1})
and changing variables from $\eta_{<j}$ to $\eta^1$, from $\eta_j$ to a bond
$b$, and from $\eta_{>j}$ to $\eta^2$, we obtain
\begin{align}\label{2nd-1st^2:eq2}
&\E[\hat\chi_{h,\beta,\Xvec}(x)^2]-\E[\hat\chi_{h,\beta,
 \Xvec}(x)]^2\nn\\
&=\sum_{\substack{\omega\in\Omega(x)\\ \eta^1\circ b\circ\eta^2\in\Omega(x)}}
 \E\bigg[e^{H_{\Xvec}(\omega)+H_{\Xvec}(\eta^1)}\E_{\Yvec}\Big[e^{-h}
 \Big(e^{-\beta X_b}-e^{-\beta Y_b}\Big)e^{H_{\Yvec}(\eta^2)}\Big]\bigg],
\end{align}
where $\eta^1\circ b\circ\eta^2$ is the concatenation of those three
paths, whose lengths are not fixed any more (due to the sum over $j$).
Since $b$ is not contained in $\eta^2$, $Y_b$ is independent of
$H_{\Yvec}(\eta^2)$, hence
\begin{align}
\E_{\Yvec}\Big[\Big(e^{-\beta X_b}-e^{-\beta Y_b}\Big)e^{H_{\Yvec}(\eta^2)}
 \Big]
&=\E_{\Yvec}\Big[e^{-\beta X_b}-e^{-\beta Y_b}\Big]~\E_{\Yvec}\Big[
 e^{H_{\Yvec}(\eta^2)}\Big]\nn\\
&= \big( e^{-\beta X_{b}}-\lambda_{\beta} \big)~\E_{\Yvec}\Big[e^{H_{\Yvec}(\eta^2)}\Big].
\end{align}
Substituting this back into (\ref{2nd-1st^2:eq2}) yields
\begin{align}\label{2nd-1st^2:eq3}
&\E[\hat\chi_{h,\beta,\Xvec}(x)^2]-\E[\hat\chi_{h,\beta,\Xvec}
 (x)]^2\nn\\[5pt]
&=e^{-h}\sum_{\substack{\omega\in\Omega(x)\\ \eta^1\circ b\circ\eta^2\in
 \Omega(x)}}\E\Big[\underbrace{e^{H_{\Xvec}(\omega)+H_{\Xvec}(\eta^1)}\big(
 e^{-\beta X_b}-\lambda_\beta\big) }_{0\text{ if }b\notin\omega}\Big]~\E_{\Yvec}
 \Big[e^{H_{\Yvec}(\eta^2)}\Big]\nn\\
&\le e^{-2h}\big(\lambda_{2\beta}-\lambda_\beta^2\big)~\E[\hat\chi_{h,\beta,
 \Xvec}(x)]\sum_{\substack{\omega^1\circ b\circ\omega^2\in\Omega(x)\\ \eta^1
 \circ b\in\Omega(x)}}\E\Big[e^{H_{\Xvec}(\omega^1)+H_{\Xvec}(\omega^2)
 +H_{\Xvec}(\eta^1)}\Big],
\end{align}
where the restricted sum over $\eta^2$ is bounded above by
$\E[\hat\chi_{h,\beta,\Xvec}(x)]$, which is translation invariant
and independent of $x\in\Zd$.

Next, we investigate the remaining sum
\begin{align}\label{2nd-1st^2:eq4}
\sum_{\substack{\omega^1\circ b\circ\omega^2\in\Omega(x)\\ \eta^1\circ
 b\in\Omega(x)}}\E\Big[e^{H_{\Xvec}(\omega^1)+H_{\Xvec}(\omega^2)
 +H_{\Xvec}(\eta^1)}\Big]\big(\ind{\omega^2\cap\eta^1=\varnothing}
 +\ind{\omega^2\cap\eta^1\ne\varnothing}\big).
\end{align}
Due to the independence among the variables in $\Xvec$,
the contribution from $\ind{\omega^2\cap\eta^1=\varnothing}$ is bounded by
\begin{align}\label{2nd-1st^2:eq5}
&\sum_{\substack{\omega^1\circ b\circ\omega^2\in\Omega(x)\\ \eta^1\circ
 b\in\Omega(x)}}\E\Big[e^{H_{\Xvec}(\omega^1)+H_{\Xvec}(\eta^1)}\Big]
 ~\E\Big[e^{H_{\Xvec}(\omega^2)}\Big]\nn\\
&\quad\le\E[\hat\chi_{h,\beta,\Xvec}(x)]\sum_{\substack{\omega^1\circ
 b\in\Omega(x)\\ \eta^1\circ b\in\Omega(x)}}\E\Big[e^{H_{\Xvec}
 (\omega^1)+H_{\Xvec}(\omega^2)}\Big]\nn\\
&\quad\le\E[\hat\chi_{h,\beta,\Xvec}(x)]\,2dB_1.
\end{align}
To bound the contribution from $\ind{\omega^2\cap\eta^1\ne\varnothing}$ in
(\ref{2nd-1st^2:eq4}), we split $\omega^2$ as $\omega^3\circ\omega^4$ at the
last visit to $\eta^1$, so that $\omega^4\cap\eta^1=\{\omega^4_0\}$.  Then, by
using the independence among the variables in $\Xvec$, we can bound the
sum over $\omega^4$ by $\E[\hat\chi_{h,\beta,\Xvec}(x)]$.  As a result, the
contribution from $\ind{\omega^2\cap\eta^1\ne\varnothing}$ is bounded by
\begin{align}\label{2nd-1st^2:eq6}
&\sum_{y\in\Zd}\sum_{\substack{\omega^1\circ b\circ\omega^3\in\Omega(x,y)\\
 \omega^4\in\Omega(y)}}\ind{\omega^1\circ b\circ\omega^3\circ\omega^4\in\Omega
 (x)}\sum_{\substack{\eta^3\in\Omega(x,y)\\ \eta^4\circ b\in\Omega(y)}}\ind{
 \eta^3\circ\eta^4\circ b\in\Omega(x)}~\ind{\omega^4\cap(\eta^3\circ\eta^4)
 =\{y\}}\nn\\
&\hskip7pc\times\E\Big[e^{H_{\Xvec}(\omega^1)+H_{\Xvec}(\omega^3)
 +H_{\Xvec}(\eta^3)+H_{\Xvec}(\eta^4)}\Big]~\E\Big[e^{H_{\Xvec}
 (\omega^4)}\Big]\nn\\[5pt]
&\le\E[\hat\chi_{h,\beta,\Xvec}(x)]\sum_{y\in\Zd}\sum_{\substack{
 \omega^1\circ b\circ\omega^3\in\Omega(x,y)\\ \eta^3\in\Omega(x,y)\\ \eta^4
 \circ b\in\Omega(y)}}\!\!\!\ind{b\notin\eta^3}\,\E\Big[e^{H_{\Xvec}(\omega^1)
 +H_{\Xvec}(\omega^3)+H_{\Xvec}(\eta^3)+H_{\Xvec}(\eta^4)}\Big]\nn\\
&=e^h\lambda_\beta^{-1}\E[\hat\chi_{h,\beta,\Xvec}(x)]\sum_{y,z\in\Zd}
 \sum_{\substack{\omega^1\in\Omega(x,z)\\ b\circ\omega^3\in\Omega(z,y)\\
 \eta^3\in\Omega(x,y)\\ \eta^4\in\Omega(y,z)}}\E\Big[e^{H_{\Xvec}
 (\omega^1)+H_{\Xvec}(b\circ\omega^3)+H_{\Xvec}(\eta^3)+H_{\Xvec}(\eta^4)}\Big]
 \nn\\
&=e^h\lambda_\beta^{-1}\E[\hat\chi_{h,\beta,\Xvec}(x)]\,B_2.
\end{align}

Finally, by summarizing (\ref{2nd-1st^2:eq3})--(\ref{2nd-1st^2:eq6}), we arrive at
\begin{align}
\frac{\E[\hat\chi_{h,\beta,\Xvec}(x)^2]-\E[\hat\chi_{h,\beta,\Xvec}(x)]^2}
 {\E[\hat\chi_{h,\beta,\Xvec}(x)]^2}\le e^{-2h}(2dB_1+e^h\lambda_\beta^{-1}B_2)
 (\underbrace{\lambda_{2\beta}-\lambda_{\beta}^2}_{O(\beta^2)}),
\end{align}
which proves (\ref{2nd-1st^2:cond}).  This completes the proof of Proposition~\ref{prop:mfbd}
\QED

\section*{Acknowledgements}
The authors are deeply indebted to two anonymous referees for their constructive
comments and numerous suggestions to improve presentation.  We would also like
to thank Rongfeng Sun for many valuable suggestions, Hubert Lacoin for clarifying
some of the details in his paper \cite{lacoin2014non} and Hugo Duminil-Copin for
pointing out typos in a previous version of the manuscript.  The first-named
author gave a talk at the IMS workshop held in Singapore during May 4--15, 2015,
and received inspiring feedback from participants.  Finally we are grateful to
Satoshi Handa and Dai Kawahara for their continual involvement in this project.

\end{document}